
\documentstyle{amsppt}
\baselineskip18pt
\magnification=\magstep1
\pagewidth{30pc}
\pageheight{45pc}
\hyphenation{co-deter-min-ant co-deter-min-ants pa-ra-met-rised
pre-print pro-pa-gat-ing pro-pa-gate
fel-low-ship Cox-et-er dis-trib-ut-ive}
\def\leaderfill{\leaders\hbox to 1em{\hss.\hss}\hfill}

\

\def\Aut{\text {\rm Aut}}

\def\sgn{{\text {\rm \, sgn}}}

\def\a{{\alpha}}
\def\be{{\beta}}
\def\g{{\gamma}}

\def\l{{\lambda}}

\def\bn{{\bold n}}

\def\bt{{\bold t}}
\def\bu{{\bold u}}
\def\bv{{\bold v}}

\def\bx{{\bold x}}

\def\b0{\text{\bf 0}}

\def\ra{{\ \longrightarrow \ }}

\def\real{{\Bbb R}}
\def\complex{{\Bbb C}}
\def\zed{{\Bbb Z}}

\def\Im{\text{\rm Im}}

\def\pd{\partial}
\def\boxit#1{\vbox{\hrule\hbox{\vrule \kern3pt
\vbox{\kern3pt\hbox{#1}\kern3pt}\kern3pt\vrule}\hrule}}
\def\rabbit{\vbox{\hbox{\kern0pt
\vbox{\kern0pt{\hbox{---}}\kern3.5pt}}}}

\def\tableau#1{
        \hbox {
                \hskip -10pt plus0pt minus0pt
                \raise\baselineskip\hbox{
                \offinterlineskip
                \hbox{#1}}
                \hskip0.25em
        }
}

\def\tabCol#1{
\hbox{\vtop{\hrule
\halign{\strut\vrule\hskip0.5em##\hskip0.5em\hfill\vrule\cr\lower0pt
\hbox\bgroup$#1$\egroup \cr}
\hrule
} } \hskip -10.5pt plus0pt minus0pt}

\def\CR{
        $\egroup\cr
        \noalign{\hrule}
        \lower0pt\hbox\bgroup$
}



\def\blank#1#2{
\hbox to #1{\hfill \vbox to #2{\vfill}}
}


\def\strut{\vrule height10pt depth5pt width0pt}

\def\hcpn{\text{\rm h}\gamma_n}
\def\hcgn{{1 \over 2} H_n}
\def\rank{\text{\rm rank}}

\def\seca{1}

\def\secb{2}
\def\secba{2.1}
\def\secbb{2.2}
\def\secbc{2.3}
\def\secc{3}
\def\seccb{3.1}
\def\seccc{3.2}
\def\seccd{3.3}
\def\secd{4}
\def\secda{4.1}
\def\secdb{4.2}
\def\secdc{4.3}

\topmatter
\title Homology representations arising from the half cube
\endtitle

\author R.M. Green \endauthor
\affil Department of Mathematics \\ University of Colorado \\
Campus Box 395 \\ Boulder, CO  80309-0395 \\ USA \\ {\it  E-mail:}
rmg\@euclid.colorado.edu \\
\newline
\endaffil


\abstract
We construct a CW decomposition $C_n$ of the $n$-dimensional half 
cube in a manner compatible with its structure as a polytope.  For each
$3 \leq k \leq n$, the complex $C_n$ has a subcomplex $C_{n, k}$, which
coincides with the clique complex of the half cube graph if $k = 4$.
The homology of $C_{n, k}$ is concentrated in degree $k-1$ and furthermore,
the $(k-1)$-st Betti number of $C_{n, k}$ is equal to the $(k-2)$-nd Betti 
number of the complement of the $k$-equal real hyperplane arrangement.
These Betti numbers, which also appear in theoretical computer science, 
numerical analysis and engineering, are the coefficients of a certain 
Pascal-like triangle (Sloane's sequence A119258).
The Coxeter groups of type $D_n$ act naturally on the complexes $C_{n, k}$,
and thus on the associated homology groups.
\endabstract

\endtopmatter

\centerline{\bf To appear in Advances in Mathematics}

\head \seca. Introduction \endhead

The half cube, also known as the demihypercube, 
may be constructed by selecting one point from each adjacent pair of 
vertices in the $n$-dimensional hypercube and taking the convex hull of 
the resulting set of $2^{n-1}$ points.  Although the resulting polytope 
is not a regular polytope like the hypercube, it still has a large symmetry
group, and its $k$-faces are of two types: regular simplices and isometric 
copies of half cubes of lower dimension.

In this paper, we will give a detailed description of the faces of the half
cube (Theorem \secbc.6) and their intersections with each other (Theorem
\secbc.8).  Using these results, we show (Theorem \seccb.2) 
that the faces of the half cube
may be assembled into a regular CW complex, $C_n$, in a natural way.
(Such a structure is sometimes called a {\it polytopal complex}.)

It is not very difficult to show that the half cube is contractible, and
therefore that
the reduced homology of the complex $C_n$ is trivial.  However,
$C_n$ has some topologically interesting subcomplexes $C_{n, k}$ for 
$3 \leq k \leq n$.  The complex $C_{n, k}$ is not a truncation of $C_n$;
rather, it is constructed from $C_n$ by deleting the interiors of all the 
half cube shaped faces of dimensions $l \geq k$.  If $k = 3$ or $k = 4$,
the complexes $C_{n, k}$ are simplicial; indeed, the case $k = 4$ gives
the clique complex of the half cube graph $\hcgn$.

Using a combination of Forman's discrete Morse theory \cite{{\bf 11}} 
and the theory of CW complexes, we prove (Theorem \seccd.2) that the 
reduced homology of $C_{n, k}$ is concentrated in degree $k-1$.  We also
give an explicit formula for the nonzero Betti numbers (Theorem \secda.2),
which also appear as the entries $T(n, n-k)$ in a certain Pascal-like triangle
(Definition \secda.3).  Although this triangle is not particularly well known,
its entries appear in a surprisingly diverse range of contexts, including:
\item{(i)}{in the problem of finding, given $n$ real numbers,
a lower bound for the complexity 
of determining whether some $k$ of them are equal \cite{{\bf 5}, {\bf 6},
{\bf 7}, \S1},}
\item{(ii)}{as the $(k-2)$-nd Betti numbers of the $k$-equal real hyperplane
arrangement in $\real^n$ \cite{{\bf 7}},}
\item{(iii)}{as the ranks of $A$-groups appearing in combinatorial homotopy
theory \cite{{\bf 1}, {\bf 2}},}
\item{(iv)}{as the number of nodes used by the Kronrod--Patterson--Smolyak
cubature formula in numerical analysis \cite{{\bf 19}, Table 3}, and}
\item{(v)}{(when $k = 3$) in engineering, as the number of three-dimensional 
block structures associated to $n$ joint systems in the construction
of stable underground structures \cite{{\bf 18}}.}

Although the relationships between (i), (ii) and (iii) above are well 
understood (see the remarks in \S\secdc\ below), the connections with the 
half cube polytope, numerical analysis and engineering appear not to have 
been noticed before.  Curiously, although explicit formulae for the numbers 
$T(n, n-k)$ are given (or implicit) in applications (ii), (iii) and (iv) 
above, these formulae all differ from each other and from our formula in
Theorem \secda.2.

In this paper, we concentrate mostly on case (ii) above.
The $k$-equal real hyperplane arrangement 
$V_{n, k}^\real$ is the set of points $(x_1, \ldots, x_n) \in \real^n$
such that $x_{i_1} = x_{i_2} = \cdots = x_{i_k}$ for some set of indices
$1 \leq i_1 < i_2 < \cdots < i_k \leq n$.  The 
complement $\real^n - V_{n, k}^\real$, denoted by $M_{n, k}^\real$, is
a manifold whose homology is concentrated in degrees $t(k-2)$, where 
$t \in \zed$ satisfies $0 \leq t \leq {\lfloor {n \over k} \rfloor}$ 
(see \cite{{\bf 7}, Theorem 1.1(b)}).  We will prove in Theorem \secda.5 and 
Corollary \secda.6 that the $(k-1)$-st Betti number of 
$C_{n, k}$ is equal to the $(k-2)$-nd Betti number of the complement of 
the $k$-equal real hyperplane arrangement.

The $n$-dimensional half cube has a large symmetry group $G_n$ of orthogonal
transformations acting on it.  This group always contains the Coxeter
group $W(D_n)$, which has order $2^{n-1}n!$, although this containment
is proper if $n = 4$.  Our final main result, Theorem \secdb.3, describes
the orbits of these groups acting on the $k$-faces of $\hcpn$; this is
useful because it induces an action of $W(D_n)$ on the nonzero homology
groups.  

The results of this paper raise various interesting questions,
which we discuss in the concluding section.

\head \secb. The geometry of the half cube \endhead

The purpose of \S\secb\ is to obtain a detailed understanding of the
geometry of the half cube, meaning a classification of its vertices, edges,
and other $k$-dimensional faces (Theorem \secbc.6).  This enables us to
verify that the faces of the half cube intersect in a nice way (Theorem
\secbc.8).  Both these results play a key role in the topological constructions
of \S\secc.  A related combinatorial problem that we solve along the way
(Proposition \secbb.6) is the classification of the cliques in the half cube
graph, which can be characterized as the $1$-skeleton of the half cube.

\subhead \secba. Polytopes \endsubhead

Following Coxeter \cite{{\bf 9}, \S7.4}, we define an $n$-dimensional 
{\it (Euclidean) polytope} 
$\Pi_n$ to be a closed, bounded, convex subset of $\real^n$ enclosed by a 
finite number of hyperplanes.  (The polytopes in this paper are all
Euclidean.)  The part of $\Pi_n$ that lies in one of
the hyperplanes is called a {\it facet}, and each facet is an 
$(n-1)$-dimensional polytope.  (Coxeter \cite{{\bf 9}} uses the term ``cell''
to mean ``facet''.)  Iterating this construction gives
rise to a set of $k$-dimensional polytopes $\Pi_k$ (called {\it $k$-faces}) 
for each $0 \leq k \leq n$.  The elements of $\Pi_0$ are called {\it vertices}
and the elements of $\Pi_1$ are called {\it edges}.  Two vertices are
called {\it adjacent} if they share an edge.  The cardinality of
the set $\Pi_i$ is denoted by $N_i$.

It is immediate from the definitions that a polytope is the convex
closure of its facets, and iteration of this observation shows that a polytope
is the convex hull of its set of vertices; in particular, the set of vertices
determines the polytope.  We can therefore speak of ``the polytope on the 
set $V$'' to refer to the polytope $\Pi(V)$ whose vertex set is $V$.

Another immediate consequence of the definitions is that the facets of a 
polytope lie in the boundary.  Conversely, the fact that there is a finite
number of bounding hyperplanes means that the boundary of a polytope 
consists precisely of the union of its facets; in other words, if a point
of the polytope lies in no facet, then it must be an interior point.

\example{Example \secba.1}
Let $x_1, \ldots, x_n$ be the usual coordinate functions in $\real^n$.
The $2n$ hyperplanes of the form $x_i = \pm 1$ for $1 \leq i \leq n$
bound a closed convex subset of $\real^n$ containing the origin.  The
corresponding polytope is the {\it hypercube} or {\it measure polytope},
denoted by $\g_n$ in Coxeter's notation \cite{{\bf 9}}.  The vertices of
$\g_n$ are the $2^n$ points of the form $
(\pm 1, \pm 1, \ldots, \pm 1)
.$
\endexample

If every $2$-dimensional face of a polytope has an even number of sides,
as is the case with $\g_n$, 
one can apply a general procedure (described in detail in 
\cite{{\bf 9}, \S8.6})
known as {\it alternation}.  This involves selecting precisely half the 
vertices of $\g_n$ in such a way that one vertex is selected from each
adjacent pair.  (For our purposes, we will include a vertex in our selection
if the entry $-1$ occurs an even number of times in its position vector.)
One can then introduce a new bounding hyperplane for each rejected vertex, and
the part of the original polytope on the same side of the hyperplane as the
rejected vertex is discarded.  This procedure does not introduce any new
vertices.

In the case of $\g_n$, the resulting polytope is called the {\it half cube}
or {\it demihypercube}, and denoted by $\hcpn$.  It is the convex hull of 
the $2^{n-1}$ vertices of $\g_n$ that contain an even number of minus signs.

\subhead \secbb. Half cube graphs \endsubhead

Let $n \geq 4$ be an integer.  We define the set $\Psi_n$ to be the set of
$2^n$ vertices of the polytope $\g_n$, $\Psi^+_n$ to be the set of
$2^{n-1}$ vertices with an even number of negative coordinates, 
and $\Psi^-_n$ to be $\Psi_n \backslash \Psi^+_n$.

The {\it Hamming distance}, $d(x, y)$, between two elements $x, y \in \Psi_n$,
is defined to be the number of coordinates at which the $n$-tuples $x$ and
$y$ differ; in other words, we have $$
2 d(x, y) = \sum_{i = 1}^n |x_i - y_i|
.$$  There is an equivalence relation $\sim$ on $\Psi_n$ given by the 
conditions that $x \sim y$ if and only if $d(x, y)$ is even, and the 
equivalence classes are precisely $\Psi^+_n$ and $\Psi^-_n$.

The {\it half cube graph}, $\hcgn$, is the simple undirected graph whose 
vertices are the elements of 
$\Psi^+$.  There is an edge between vertices $x$ and $y$ if and only if 
$d(x, y) = 2$.  (Note that we could equally well have defined this graph using
$\Psi^-$.)

A {\it $k$-clique} (or {\it clique} for short) in a graph $G$ is a set of $k$
vertices of $G$, any two of which are mutually adjacent.  The {\it clique
complex} of $G$ is the abstract simplicial complex whose $k$-faces are the
$k$-cliques of $G$.  

\definition{Definition \secbb.1}
Let $\bn = \{1, 2, \ldots, n\}$, $\bv \in \Psi^-$ and $S \subseteq \bn$.
If $|S| = k > 0$, we define a $k$-clique $K(\bv, S)$ of $\hcgn$ by the 
condition that $\bv' \in K(\bv, S)$ if and only if $\bv$ and $\bv'$ differ
only in the $i$-th coordinate, and $i \in S$.  (It is immediate from the
definitions that $K(\bv, S)$ is a $k$-clique.)  We call $S$ a {\it mask}
of $K(\bv, S)$.  We say that $K(\bv, S)$ is {\it opposite $\bv$}, or that
$\bv$ is an {\it opposite point} of $K(\bv, S)$. 
\enddefinition

\proclaim{Lemma \secbb.2}
Let $K(\bv, S)$ be a set of the form given in Definition \secbb.1.  If
$|K(\bv, S)| \geq 3$, then the point $\bv$ and the set $S$ are each 
determined by the set $K(\bv, S)$.
\endproclaim

\demo{Proof}
Suppose $K$ is a set of the required form with $|K| \geq 3$.  For each
$1 \leq i \leq n$, it is the case that a majority of the points in $K$
(all but one if $i \in S$, and all of them if $i \not\in S$) contain 
the entry $a_i$ in the $i$-th coordinate
for some $a_i \in \{\pm 1\}$.  We define $$
\bv = (a_1, a_2, \ldots, a_n)
.$$  The subset $S \subseteq \bn$ is then the set of all $j$ for which
some point of $K$ differs from $\bv$ only in the $j$-th coordinate.
\qed\enddemo

\example{Example \secbb.3}
Let $n = 5$ and $\bv_1 = (1, -1, -1, 1, -1)$.  Then we have $$
\eqalign{
K(\bv_1, \{3\}) &= \{(1, -1, 1, 1, -1)\},\cr
K(\bv_1, \{3, 4\}) &= \{(1, -1, 1, 1, -1), (1, -1, -1, -1, -1)\},\cr
K(\bv_1, \{2, 3, 4\}) &= 
\{(1, 1, -1, 1, -1), (1, -1, 1, 1, -1), (1, -1, -1, -1, -1)\}.\cr
}$$  In contrast to Lemma \secbb.2, we have $K(\bv_1, \{3, 4\}) 
= K(\bv_2, \{3, 4\})$, where $\bv_2 = (1, -1, 1, -1, -1)$, so that 
$K(\bv_1, \{3, 4\})$ has two
opposite points, but only one mask.  There are a total of $5$ pairs, 
$(\bv_1, S)$ giving rise to the singleton set $K(\bv_1, \{3\})$, and the set 
has $5$ opposite points, each with a different mask.
\endexample

\definition{Definition \secbb.4}
Let $\bn = \{1, 2, \ldots, n\}$, $\bv \in \Psi^+$ and let $S \subseteq \bn$.
We define the subset $L(\bv, S)$ of $\hcgn$ by the 
condition that $\bv' \in L(\bv, S)$ if and only if $\bv$ and $\bv'$ agree
in the $i$-th coordinate whenever $i \not\in S$.  The set $S$ is called the
{\it mask} of $L(\bv, S)$; it can be characterized as the set of coordinates
at which not all points of $L(\bv, S)$ agree.

In particular, if $|S| = 3$, then the elements of $L(\bv, S)$
are $\bv$ itself, together with the three vectors that differ from
$\bv$ in precisely two coordinates indexed by $S$, and it follows from the
definitions that $L(\bv, S)$ is a $4$-clique.  
\enddefinition

\example{Example \secbb.5}
Let $n = 5$, $\bv = (1, -1, -1, 1, 1)$ and $S = \{1, 3, 4\}$.  Then we have $$
L(\bv, S) = \{(1, -1, -1, 1, 1), (-1, -1, 1, 1, 1), (-1, -1, -1, -1, 1),
(1, -1, 1, -1, 1)\}
.$$
\endexample

Note that a $4$-clique of the form $L(\bv, S)$ cannot be of the form 
$K(\bv', S')$ for any $\bv'$ and $S'$: the elements of a $4$-clique of 
the form $K(\bv', S')$ (respectively, $L(\bv, S)$) agree at all but precisely
four (respectively, three) coordinates.

\proclaim{Proposition \secbb.6}
Every clique in the graph $\hcgn$ is of the form $K(\bv, S)$ for some $S$,
or of the form $L(\bv, S)$ for some $S$ with $|S| = 3$.
\endproclaim

\demo{Proof}
It is enough to show that if $C$ is a $k$-clique of one of the above forms,
then any extension of $C$ to a $(k+1)$-clique $C'$ is also of one of the 
above forms.

If $k = 1$, then $C$ is necessarily a single point $\bv_1$, and $C'$
consists of two points $\bv_1$ and $\bv_2$ which differ from each other in
precisely two coordinates.  This satisfies the required condition.

If $k = 2$, then $C = \{\bv_1, \bv_2\}$ consists of two points differing
from each other in precisely two coordinates, $a$ and $b$.  Suppose that
$C' = \{\bv_1, \bv_2, \bv_3\}$ is an extension of $C$ to a $3$-clique.
Since $\{\bv_1, \bv_3\}$ is a $2$-clique, it must be the case that $\bv_3$ 
differs from $\bv_1$ in precisely two coordinates, $c$ and $d$.  However,
since $\{\bv_2, \bv_3\}$ is a $2$-clique, it follows that $\bv_2$ and
$\bv_3$ differ from each other in precisely two coordinates, and thus that
the set $\{a, b, c, d\}$ has cardinality $3$.  Without loss of generality,
we may assume that $d = a$ and $b \ne c$.  Letting $\bv' \in \Psi^-_n$ be
the unique element that differs from $\bv_1$ only in the $a$-th coordinate,
we find that $C' = K(\bv', \{a, b, c\})$, which satisfies the required
conditions.

If $k = 3$, then the preceding arguments show that 
$C = \{\bv_1, \bv_2, \bv_3\}$ is of the form
$K(\bv', \{a, b, c\})$, where $a$, $b$ and $c$ are as in the previous
paragraph.  Suppose that $C' = \{\bv_1, \ldots, \bv_4\}$ is 
a $4$-clique extending $C$.  If $\bv_4$ differs from $\bv_1$ in the 
$a$-th coordinate, then it also differs from $\bv_1$ in the $d$-th coordinate
for some $d \not\in \{a, b, c\}$, and we have $C' = K(\bv', \{a, b, c, d\})$.
The other possibility is that $\bv_4$ agrees with $\bv_1$ in the $a$-th
coordinate, which means that it disagrees with each of $\bv_2$ and $\bv_3$
in the $a$-th coordinate.  
Because $\bv_4 \ne \bv_1$ and $\{\bv_2, \bv_4\}$ is a $2$-clique, it must be 
the case that $\bv_4$ agrees with $\bv_2$ in the $b$-th coordinate, and
a similar argument shows that $\bv_4$ agrees with $\bv_3$ in the $c$-th
coordinate.  The fact that $\{\bv_1, \bv_4\}$ is a $2$-clique then shows
that $\bv_4$ agrees with $\bv_1$ at all but two coordinates: $b$ and $c$.
In summary, we have $C' = L(\bv_1, \{a, b, c\})$, completing the analysis of
the case $k = 3$.

Suppose that $k = 4$.  We will show that, in the above notation, if \newline
$C = \{\bv_1, \ldots, \bv_4\}$ is of the form $L(\bv_1, \{a, b, c\})$, then 
$C$ cannot be extended to a $5$-clique $C' = \{\bv_1, \ldots, \bv_5\}$.
If this were the case, then the $4$-clique $\{\bv_1, \bv_2, \bv_3, \bv_5\}$ 
would be an extension of $\{\bv_1, \bv_2, \bv_3\}$, and the argument of
the previous paragraph combined with the fact that $\bv_4 \ne \bv_5$ then
shows that $\bv_5$ would differ from $\bv_1$ in two coordinates, $a$ and $e$,
where $e \not\in \{a, b, c\}$.  This would imply that $\bv_4$ and $\bv_5$
would differ in four coordinates, $a$, $b$, $c$ and $e$, contradicting the
requirement that $\{\bv_4, \bv_5\}$ be a $2$-clique.  

We may therefore assume that $C = \{\bv_1, \ldots, \bv_4\}$ is of the form 
$K(\bv', \{a, b, c, d\})$, and that each of the points $\bv_2, \bv_3, \bv_4$
differs from $\bv_1$ in two coordinates, one of which is $a$.  Suppose again 
that the $5$-clique  $C' = \{\bv_1, \ldots, \bv_5\}$ is an extension of $C$.
Since $\{\bv_1, \bv_2, \bv_3, \bv_5\}$ is a clique of the form 
$K(\bv'', S)$, the argument used to deal with the case $k = 3$ shows that
$\bv_5$ differs from $\bv_1$ in two coordinates, one of which is $a$.  We
define $e$ to be the other coordinate at which these points differ; since
$C'$ has cardinality $5$, it must be the case that the set $S = \{a, b, c, d,
e\}$ also has cardinality $5$.  This proves that $C' = K(\bv', S)$, and
completes the proof of the case $k = 4$.

The case of $k > 4$ proceeds similarly.  Let $C = \{\bv_1, \ldots, \bv_k\}$
be a clique of size $k > 4$; by the arguments above, it follows that
$C = K(\bv', S)$ for suitable $\bv'$ and $S$.  Moreover, we can arrange
things so that each of $\bv_2$, $\bv_3$ and $\bv_4$ differ from $\bv_1$ in
the $a$-th coordinate and one other coordinate (the $b$-th, $c$-th and $d$-th,
respectively), and that $\bv'$ differs from $\bv_1$ only in the $a$-th
coordinate.  Suppose that $C' = \{\bv_1, \bv_{k+1}\}$ is an extension of 
$C$ to a $(k+1)$-clique, and consider the $4$-subclique $C'' = \{\bv_1, \bv_2,
\bv_3, \bv_l\}$ for some $4 < l \leq k+1$.  Since $C''$ must be of the form 
$K(\bv'', S')$, it must be the case that $\bv_l$ differs from $\bv_1$ in 
two coordinates, one of which is the $a$-th.  Since this is true for all
$4 < l \leq k+1$, it follows that $C' = K(\bv', S'')$, where $\bv'$ is
as before, and $S''$ consists of all coordinates at which some pair of
elements in $C'$ disagree.  This completes the proof.
\qed\enddemo

\remark{Remark \secbb.7}
The above proof shows that, if $n > 4$, a $4$-clique is of the form 
$K(\bv, S)$ if it can be extended to a $5$-clique, and is of the form
$L(\bv, S)$ otherwise.  In particular, there can be no automorphism of
the graph $\hcgn$ exchanging $4$-cliques of different types.
\endremark

\subhead \secbc. Faces of the half cube \endsubhead

In \S\secbc, we give a more explicit description of the half cube $\hcpn$.
The vertices of this polytope were described in \S\secba\  above, and we
will now describe the $k$-faces for $0 < k < n$.  Although such a description
is not new, the only reference we know for the result of Theorem \secbc.6
is the proof of \cite{{\bf 10}, Proposition 4.2}, and unfortunately, the latter 
paper does not give a detailed enough statement of the result for our 
purposes, nor does it give the details of the proof.

Recall in the sequel that a polytope is the convex hull of its set
of vertices.

\proclaim{Lemma \secbc.1}
\item{\rm (i)}
{Let $\bv' = (v_1, \ldots, v_n) \in \Psi^-_n$ and $S \subseteq \bn$.  The 
convex hull of the set $K(\bv', S)$ consists precisely of the points 
$(x_1, \ldots, x_n) \in \real^n$ satisfying the following three conditions:
\item{\rm (a)}{$x_i = v_i$ for all $i \not\in S$;}
\item{\rm (b)}{$\sgn(v_i) (x_i - v_i) \leq 0$ for all $1 \leq i \leq n$,
where $\sgn(v_i) = v_i/|v_i|$; and}
\item{\rm (c)}{$\sum_{i \in S} \sgn(v_i) (x_i - v_i) = -2$.}}
\item{\rm (ii)}
{Let $\bx = (x_1, \ldots, x_n)$ be a point in the convex hull of the set 
$L(\bv, S)$ for some $\bv \in \Psi^+_n$ and $S \subseteq \bn$.  If $x_i = c$
for $c \in \{1, -1\}$ then $\bx$ is a convex combination of a subset of 
points of $L(\bv, S)$ all of which agree with $\bx$ in the $i$-th coordinate.}
\endproclaim

\demo{Proof}
The convex hull of the set $K(\bv', S)$ is the set of convex linear
combinations of $K(\bv', S)$, that is, the set of all points of the form $$
\sum_{\bu \in K(\bv', S)} \l_\bu \bu
$$ such that $\l_\bu \geq 0$ for all $\bu$ and such that $$
\sum_{\bu \in K(\bv', S)} \l_\bu = 1
.$$  It is convenient to shift the origin of the problem by subtracting
$\bv'$ from all points.  The problem then becomes to determine the convex
hull of the points $$
\{-2 \sgn(v_i) e_i : i \in S
\},$$ where $\{e_1, \ldots, e_n\}$ is the standard basis for $\real^n$.
This second convex hull is equal to $$
\left\{ (x_1, \ldots, x_n) \in \real^n : 
x_i = 0 \text{\rm \ for all\ } i \not\in S
\text{\rm \ and\ } 
\sum_{i \in S} {{-1} \over {2 \sgn(v_i)}} x_i = 1 \right\}
,$$ where each term in the sum is nonnegative.  An equivalent description of
this set is $$
\left\{ (x_1, \ldots, x_n) \in \real^n : 
x_i = 0 \text{\rm \ for all\ } i \not\in S
\text{\rm \ and\ } 
\sum_{i \in S} \sgn(v_i) x_i = -2 \right\}
,$$ where $\sgn(v_i) x_i \leq 0$ for all $i \in S$, and assertion (i) follows.

For (ii), consider a convex combination $$
\bx = \sum_{\bu \in L(\bv, S)} \l_\bu \bu
$$ such that $\l_\bu \geq 0$ for all $\bu$ and such that $$
\sum_{\bu \in L(\bv, S)} \l_\bu = 1
.$$  The fact that there are only two possible entries in the $i$-th 
coordinate for each $\bu$, namely $\pm 1$, means that the only way we can
have $x_i = c$ is for all $\bu = (u_1, \ldots, u_n)$ such that 
$\l_\bu \ne 0$ to satisfy $u_i = c$.  The assertion of (ii) follows from this.
\qed\enddemo

\proclaim{Lemma \secbc.2}
\item{\rm (i)}
{If $n = 4$, the half cube $\hcpn$ has precisely $16$ facets: eight of these
are the polytopes $\Pi(K(\bv', \{1, 2, 3, 4\}))$ arising with $4$-cliques
(where $\bv' \in \Psi^-_n$)
and the other eight are the polytopes $L(\bv, S)$ associated with the 
$4$-cliques obtained as $S$ ranges over the eight subsets of $\bn$
and $\bv \in \Psi_n^+$.}
\item{\rm (ii)}
{Suppose that $n > 4$.  Then the number of facets of the half cube 
$\hcpn$ is $2^{n-1} + 2n$.  Of these, $2^{n-1}$ facets are the polytopes
$\Pi(K(\bv', \{1, \ldots, n\}))$ arising from $n$-cliques (where $\bv'
\in \Psi^-_n$) and the other
$2n$ are the polytopes $\Pi(L(\bv, S))$, where $\bv \in \Psi^+_n$ and
$S \subset \bn$ satisfies $|S| = n-1$.}
\endproclaim

\demo{Proof}
The proof of part (i) can be regarded as special case of the proof of 
part (ii), if we define $\hcpn$
for $n = 3$ in the obvious way as the convex hull of alternating points on a
regular cube (namely $L(\bv, S)$ for suitable $\bv$ and $S$).  We will 
therefore only prove (ii).

Coxeter \cite{{\bf 9}, \S8.6} proves that the half cube $\hcpn$ has $2^{n-1}$
facets that are regular simplices, and $2n$ facets that are half cubes of
smaller dimension.  The alternation construction described in \S\secba\ 
creates no new vertices and shows that there are no other facets.  General 
properties of facets show that the vertices of facets must be vertices of the 
original polytope.  This reduces the problem to finding $2^{n-1}$ different
regular simplices and $2n$ different half cubes, all of dimension $n-1$, whose
vertices are contained in the set of vertices for $\hcpn$.

For the simplices, note that the cliques $K(\bv', \{1, \ldots, n\})$ as
$\bv'$ ranges over the points of $\Psi^-_n$ are distinct by Lemma \secbb.2.
They exhaust the possible simplices on the vertex set of $\hcpn$ by 
Proposition 2.2.6, so these must be the $2^{n-1}$ faces described in
\cite{{\bf 9}, \S8.6}.

For the half cubes, note that there are $2n$ distinct sets of the form
$L(\bv, S)$ as $S$ ranges over the subsets of $\bn$ of size $n-1$, and each of
these $2n$ sets contains $2^{n-2}$ points of $\Psi^+_n$.  (Note that for each
set, there are $2^{n-1}$ possible choices for $\bv$.)  If one deletes
the coordinate of the points of $L(\bv, S)$ that is not indexed by $S$,
then one obtains either $\Psi^+_{n-1}$ or $\Psi^-_{n-1}$, which shows that
$L(\bv, S)$ is isometric (by applying a reflection in one coordinate if
necessary) to a half cube of dimension $n-1$.  The mask $S$ of $L(\bv, S)$
is characterized as the set of coordinates at which not all $2^{n-2}$ elements
of $L(\bv, S)$ agree, from which it follows that there are $2n$ distinct
subsets of the form $L(\bv, S)$.  However, the sets $L(\bv, S)$ can also
be characterized as those vertices that lie in some fixed coordinate plane
of the form $x_i = \pm 1$, and these are bounding faces of $\hcpn$, so by
Lemma \secbc.1 (ii), the 
convex hulls of the $2n$ distinct subsets are faces of $\hcpn$ that are
also half cubes of dimension $n-1$.  They are therefore equal to the $2n$
faces described in \cite{{\bf 9}, \S8.6}, and the statement follows.
\qed\enddemo

\proclaim{Lemma \secbc.3}
Let $V' = K(\bu', \bn)$, where $\bu' \in \Psi^-_n$, so that $\Pi(V')$ is 
a facet of $\hcpn$.
Let $V$ be a set of vertices of the form $K(\bv', S)$,
where $\bv' \in \Psi^-_n$ and $S \subseteq \bn$.  
\item{\rm (i)}{The intersection $\Pi(V) \cap \Pi(V')$ is either empty, or
is another polytope of the form $\Pi(K(\bx', S'))$, where $\bx'
\in \Psi^-_n$ and $S' \subseteq \bn$.  In either case, we must 
have $K(\bx', S') = K(\bv', S) \cap K(\bu', \bn)$.}
\item{\rm (ii)}{If $V \ne V'$ and some interior point of $\Pi(V)$ lies in 
$\Pi(V')$, then it must be the case that $V \subset V'$ and that $\Pi(V)$
is contained in the boundary of $\Pi(V')$.}
\endproclaim

\demo{Proof}
Suppose that $\bx = (x_1, \ldots, x_n)$ is a point in the intersection of
the two polytopes.  Since $\bv'$ and $\bu'$ lie in $\Psi^-_n$, they
must differ in an even number of coordinates.
If $\bv' = \bu'$, then $K(\bv', S) \subseteq K(\bu', \bn)$ and
(i) and (ii) follow.  If $\bv' \ne \bu'$, let $i$ and $j$ be two positions
at which $\bv'$ and $\bu'$ differ.  In the
notation of Lemma \secbc.1 (i) with $\bu' = (u_1, \ldots, u_n)$, this means
that $\sgn(v_i) = -\sgn(u_i)$ and $\sgn(v_j) = -\sgn(u_j)$.  It follows from
this that we have $$
\sgn(v_i)(x_i - v_i) + \sgn(u_i)(x_i - u_i) = -2
$$ and $$
\sgn(v_j)(x_j - v_j) + \sgn(u_j)(x_j - u_j) = -2
.$$  Lemma \secbc.1 (i) shows that all four terms being added are nonpositive,
and the only way this can be compatible with the inequalities $$
\sum_{1 \leq i \leq n} \sgn(v_i) (x_i - v_i) = -2
$$ and $$
\sum_{1 \leq i \leq n} \sgn(u_i) (x_i - u_i) = -2
$$ and the sign constraints mentioned in Lemma \secbc.1 (i) is to have both $$
\sgn(v_i)(x_i - v_i) + \sgn(v_j)(x_j - v_j) = -2 \eqno(1)
$$ and $$
\sgn(u_i)(x_i - u_i) + \sgn(u_j)(x_j - u_j) = -2. \eqno(2)
$$  

In order for $\bx$ to exist, that is, for (1) and (2) to have any common 
solutions, we need either $i \in S$ or $j \in S$ or both.  If 
$S \cap \{i, j\} = \{i\}$, the only solution is $x_j = v_j$ and $x_i = u_i$, 
meaning that $\Pi(V) \cap \Pi(V') = \bx$, where $\bx$ is the common point
of $V$ and $V'$ that differs from $\bv'$ only in the $i$-th coordinate,
and differs from $\bu'$ only in the $j$-th coordinate.  This satisfies (i),
and (ii) holds vacuously because $\bx$ is not an interior point of $\Pi(V)$.
A similar argument deals with the case
$S \cap \{i, j\} = \{j\}$.

We may now assume that $\{i, j\} \subseteq S$, in which case
(1) and (2) force $x_k = u_k = j_k$ for all $1 \leq k \leq n$ other
than $k = i$ and $k = j$.  Let us define $\bx_1$ (respectively, $\bx_2$) to 
be the element of $\Psi^+_n$ that agrees with $\bv'$ except in the 
$i$-th (respectively, $j$-th) coordinate, so that $\bx_1$ and $\bx_2$ lie
in $K(\bv', S) \cap K(\bu', \bn)$.  We have shown that $\bx_1$ and $\bx_2$
are the only two vertices that lie in both polytopes, and it follows 
that $\bx$ lies 
on the edge containing $\bx_1$ and $\bx_2$.  This edge is therefore the
intersection of the two polytopes; we thus have $\Pi(V) \cap \Pi(V') = 
\Pi(V \cap V')$, satisfying (i).  Either this edge is the whole
of $\Pi(V)$, or it lies entirely within the boundary of $\Pi(V)$, and
in either case, (ii) holds.
\qed\enddemo

\proclaim{Lemma \secbc.4}
Let $\Pi(V')$ be a facet of $\hcpn$ with vertex set $V'$.  Let $V$ be 
a set of vertices of the form $K(\bv', S)$ or $L(\bv, T)$, where $\bv' 
\in \Psi^-_n$, $\bv \in \Psi^+_n$, $S \subseteq \bn$ and $T \subset 
\bn$ satisfies $2 \leq |T| < n$.
\item{\rm (i)}{The intersection $\Pi(V) \cap \Pi(V')$ is either empty, or
is another polytope of the form $\Pi(V'')$, where $V''$ is 
of the form $K(\bt', S')$ or $L(\bt, T')$, and we have
$\bt' \in \Psi^-_n$, $\bt \in \Psi^+_n$, $S' \subseteq \bn$ and $T \subset
\bn$.  In either case, we must have $V'' = V \cap V'$.}
\item{\rm (ii)}{If $V \ne V'$ and some interior point of $\Pi(V)$ lies in 
$\Pi(V')$, then it must be the case that $V \subset V'$ and that $\Pi(V)$
is contained in the boundary of $\Pi(V')$.}
\endproclaim

\demo{Note}
Note that the inclusion $\Pi(V \cap V') \subseteq \Pi(V) \cap \Pi(V')$ in 
(i) is immediate, but that the reverse inclusion is not.
\enddemo

\demo{Proof}
By Lemma \secbc.2, any facet of $\hcpn$ is of the form 
$\Pi(K(\bu', \bn))$ or of the form $\Pi(L(\bu, T'))$, where 
$\bu' \in \Psi^-_n$, $\bu \in \Psi^+_n$ and $T' \subset \bn$ satisfies
$|T'| = n-1$.

There are four cases to consider.  The first, and most complicated case
is an intersection of the form $$
\Pi(K(\bv', S)) \cap \Pi(K(\bu', \bn))
,$$ and this follows from Lemma \secbc.3.

The second of the four cases is an intersection of the form $$
\Pi(L(\bv, T)) \cap \Pi(K(\bu', \bn))
.$$  Suppose that $\bx$ is an element of the intersection.  If $i$ is a 
coordinate not indexed by $T$, all the points in $\Pi(L(\bv, T))$ must agree 
with $\bv$ in the $i$-th coordinate, and the same is true for $\bx$.

If $\bu'$ disagrees with $\bv$ in the $i$-th coordinate, the only solution
for $\bx$ compatible with Lemma \secbc.1 (i) is for $\bx$ to be the point
of $K(\bu', \bn)$ that disagrees with $\bu'$ only in the $i$-th coordinate.
Even this can only happen if $\bx$ agrees with $\bv$ in all other coordinates
not indexed by $T$.  This proves that $\bx$ is the unique element of
$V \cap V'$ in this case, and (i) and (ii) follow because $\bx$ cannot be
an interior point of either polytope.

We have now reduced to the case where $\bu'$ agrees with $\bv$ at all
coordinates not indexed by $T$.  In order for $\bx$ to exist, we also need
$\bx$ to agree with $\bv$ at all coordinates not indexed by $T$.
Since $$K(\bu', T) = L(\bv, T) \cap K(\bu', \bn),$$
Lemma \secbc.1 (i) and the fact that $\bx \in \Pi(K(\bu', \bn))$
now show that $\bx \in \Pi(K(\bu', T))$, which proves (i).
Lemma \secbc.2 shows
that $\Pi(K(\bu', T))$ lies in the boundary of $L(\bv, T)$, from which
assertion (ii) follows.

The third case involves an intersection of the form $$
\Pi(K(\bv', S)) \cap \Pi(L(\bu, T'))
,$$ which is similar to, but easier than, the second case.  
If $\bv'$ disagrees with
$\bu$ at the unique coordinate not indexed by $T'$, then the intersection
of the two polytopes is the unique point of intersection of the sets
$V$ and $V'$, which proves (i) and (ii) because this point cannot
be interior to $\Pi(V)$.

Suppose on the other hand that $\bv'$ agrees with $\bu$ at the 
coordinate not indexed by $T'$, and let $\bx$ be a point in the intersection
of the two polytopes.  Since $\bx$ lies in $\Pi(L(\bu, T'))$, $\bx$ agrees
with $\bu$ at the coordinate not indexed by $T'$, and it follows that
$\bx$ agrees with $\bv'$ at the coordinate not indexed by $T'$.  
Lemma \secbc.1 (i) then shows that $\bx$ lies in $\Pi(K(\bv', S \cap T'))$.
Assertion (i) follows from the observation that $$
K(\bv', S) \cap L(\bu, T') = K(\bv', S \cap T')
,$$ and assertion (ii) follows from the fact that either $S \subseteq T'$
or $\Pi(K(\bv', S \cap T'))$ lies in the boundary of the simplex 
$\Pi(K(\bv', S))$.

The fourth and final case deals with the intersection $$
\Pi(L(\bv, T)) \cap \Pi(L(\bu, T'))
;$$ recall that $|T'| = n-1$ in this case.
If $i \not\in T \cup T'$ and $\bv$ and $\bu$ disagree at the $i$-th 
position, then the intersection is empty and we are done, so suppose that 
$\bx$ is a point in the intersection.

Suppose first that $T \not\subseteq T'$.
In this case, $\bx$ must agree with each of $\bv$ at all coordinates
outside $T'$, and $\bx$ must agree with $\bu$ at all coordinates outside
$T$.  Let $\bx'$ be a point of $\Psi_n^+$ that agrees with $\bv$ (respectively,
$\bu$) at all points outside $T'$ (respectively, $T$); such a point exists
because $|T| \geq 2$ by assumption, and thus $|T \cap T'| \geq 1$.
Lemma \secbc.1 (ii) now 
shows that $\bx$ lies in $\Pi(L(\bx', T \cap T'))$.  Since $$
L(\bv, T) \cap L(\bu, T') = L(\bx', T \cap T')
,$$ we have shown that $$
\Pi(L(\bv, T)) \cap \Pi(L(\bu, T')) = \Pi(L(\bx', T \cap T'))
,$$ which proves (i) in this case.  If $i$ is the unique coordinate in
$T \backslash (T \cap T')$, then the $i$-th coordinate of $\bx$ is $\pm 1$.  
This means that $\bx$ is not an interior point of $\Pi(V)$, and (ii) follows.

The other possibility is that $T \subseteq T'$, and that $\bv$ and $\bu$
agree at all coordinates not indexed by $T'$.  In this case, we have $$
L(\bv, T) \subseteq L(\bv, T') = L(\bu, T')
$$ and $$
\Pi(L(\bv, T)) \cap \Pi(L(\bu, T')) = \Pi(L(\bv, T))
,$$ which proves (i).  Assertion (ii) now follows from
Lemma \secbc.2.
\qed\enddemo

\proclaim{Lemma \secbc.5}
Any subset of $\Psi^+_n$ of the form $K(\bv', S)$ or $L(\bv, T)$, where
$\bv' \in \Psi^-_n$, $\bv \in \Psi^+_n$, $S \subseteq \bn$ and 
$T \subsetneq \bn$, is equal to an intersection of subsets of
$\Psi^+_n$ of the form $K(\bu', \bn)$ and/or $L(\bu, T')$, where
$\bu' \in \Psi^-_n$, $\bu \in \Psi^+_n$ and $|T'| = n-1$.
\endproclaim

\demo{Proof}
Consider a subset of $\Psi$ of the form $K(\bv', S)$, where $\bv' \in \Psi^-_n$
and $S \subseteq \bn$, and observe that $K(\bv', S)$ is a subset of
$K(\bv', \bn)$.  Suppose that there exists $i \in \bn \backslash S$.  
Let $\bu_i \in \Psi^+$ be any
element agreeing with $\bv'$ in the $i$-th coordinate (note that such a $\bu$
does exist), and let $S_i = \bn \backslash \{i\}$.  It follows that $$
K(\bv', S) = K(\bv', \bn) \cap \bigcap_{i \in \bn \backslash S} L(\bu_i, S_i)
.$$  The assertion for the sets $K(\bv', S)$ follows by intersecting
the simplex shaped facet corresponding to $K(\bv', \bn)$ with the bounding
hyperplanes corresponding to the sets $L(\bu_i, S_i)$.

The assertion for the sets $L(\bv, T)$ follows similarly from the fact 
that $$
L(\bv, T) = \bigcap_{i \in \bn \backslash T} L(\bv, S_i)
,$$ where $S_i$ is as above.
\qed\enddemo

\proclaim{Theorem \secbc.6}
The $k$-faces of $\hcpn$ for $k < n$ are as follows:
\item{\rm (i)}{$2^{n-1}$ $0$-faces (vertices) given by the elements of
$\Psi^+_n$;}
\item{\rm (ii)}{$2^{n-2} {n \choose 2}$ $1$-faces $\Pi(K(\bv', S))$,
where $\bv' \in \Psi^-_n$ and $|S| = 2$;}
\item{\rm (iii)}{$2^{n-1} {n \choose 3}$ simplex shaped $2$-faces 
$\Pi(K(\bv', S))$, where $\bv' \in \Psi^-_n$ and $|S| = 3$;}
\item{\rm (iv)}{$2^{n-1} {n \choose {k+1}}$ simplex shaped $k$-faces 
$\Pi(K(\bv', S))$, where $\bv' \in \Psi^-$ and $|S| = k+1$ for $3 \leq k < n$;}
\item{\rm (v)}{$2^{n-k} {n \choose k}$ half cube shaped $k$-faces 
$\Pi(L(\bv, S))$, where $\bv \in \Psi^+_n$ and $|S| = k$ for $3 \leq k < n$.}
\endproclaim

\demo{Proof}
The numbers in (i), (ii) and (iii) are the numbers of $1$-cliques, $2$-cliques
and $3$-cliques in $\hcgn$, which correspond to simplex-shaped faces.

Suppose that $k \geq 3$.  There are $2^{n-1} {n \choose {k+1}}$ 
$(k+1)$-cliques of the form $K(\bv, S)$ in $\hcgn$: they are determined by 
Lemma \secbb.2 by their ${n \choose {k+1}}$ possible masks and by their 
$2^{n-1}$ possible opposite points.  These are the simplex-shaped faces.  
The half cube shaped faces of dimension $k$ have ${n \choose k}$ possible
masks.  For each mask there are $2^{n-k}$ values for coordinates outside
the mask.  This proves (iv) and (v) and completes the proof.
\qed\enddemo

\remark{Remark \secbc.7}
If $k = 3$, parts (iv) and (v) of Theorem \secbc.6 give 
two types of simplex shaped
$3$-faces, corresponding to the two types of $4$-cliques in $\hcpn$.  If
$n > 4$, they can be distinguished by the fact that the first type of 
$4$-clique appears as the face of a $5$-clique, and the other does not.
\endremark

\proclaim{Theorem \secbc.8}
\item{\rm (i)}
{If $\Pi(V)$ and $\Pi(V')$ are faces of $\hcpn$, where $V$ and $V'$ are as in
Theorem \secbc.6, then $\Pi(V) \cap \Pi(V') = \Pi(V \cap V')$.}
\item{\rm (ii)}
{No point of $\hcpn$ can be an interior point of more than one face.}
\item{\rm (iii)}
{Every point of $\hcpn$ is either an interior point of $\hcpn$, or a vertex
of $\hcpn$, or is an interior point of a (unique) $k$-face of $\hcpn$ for 
some $0 < k < n$.}
\endproclaim

\demo{Proof}
We first prove (i).  If $V'$ is the whole of $\Psi^+_n$, then the claim
follows immediately, so we may suppose that this is not the case.

Theorem \secbc.6 shows that $V'$ is of the
form $K(\bv', S)$ or $L(\bv, T)$, where $\bv' \in \Psi^-_n$, $\bv \in
\Psi^+_n$, $S \subseteq \bn$ and $T \subsetneq \bn$.
Lemma \secbc.5 then shows that $V'$ can be written
as a finite intersection $$
V' = \bigcap_{i = 1}^r V_i
$$ for some $r \geq 1$, where each $V_i$ is of
the form $K(\bu', \bn)$ or $L(\bu, T')$ for $\bu' \in \Psi^-_n$, 
$\bu \in \Psi^+_n$ and $|T'| = n-1$.  The proof proceeds by induction on
$r$, and the base case, $r = 1$, is Lemma \secbc.4 (i).  The inductive step
follows by writing $$
V \cap V' = \left(V \cap \bigcap_{i = 1}^{r-1} V_i \right) \cap V_r
$$ and appealing to the base case.  This completes the proof of (i).

We now turn to (ii).  Suppose for a contradiction that $\bx$ is a common
interior point of the distinct faces $\Pi(V)$ and $\Pi(V')$ (not necessarily
of the same dimension).  

If $\Pi(V')$ is the whole of $\hcpn$ then it
follows that $\Pi(V)$ is not, and hence that $\Pi(V)$ is contained in a face
of $\hcpn$.  This is a contradiction because all points of $\Pi(V)$ are
contained in the boundary of $\hcpn$, and thus cannot be interior points of
$\Pi(V')$.

We may therefore assume that $\Pi(V')$ is not the whole of $\hcpn$, and we
may write $V' = \bigcap_{i = 1}^r V_i$ as in the proof of (i), so that 
$r \geq 1$.  By part (i), we have $$
\Pi(V') = \bigcap_{i = 1}^r \Pi(V_i)
,$$ where the $\Pi(V_i)$ are facets.  Since, for each $i$, $\bx$ is an interior
point of $\Pi(V)$ that lies in $\Pi(V_i)$, Lemma \secbc.4 (ii) shows that
$V \subseteq V_i$ for all $i$, and hence that $V \subseteq V'$.  Reversing
the roles of $V$ and $V'$ in the argument then shows that $V = V'$, which
is a contradiction establishing (ii).

To prove (iii), we first remark that the analogous result for simplices
is well known (see \cite{{\bf 13}, Definition 6.2 ($2'$)}).  The uniqueness
part of (iii) is immediate from part (ii), so we concentrate on establishing
existence, and the proof is by induction on the dimension of the largest
face containing the point $\bx$ in question 
(where we consider $\hcpn$ to be an $n$-face).  The base case, involving
a polytope of dimension $0$, is trivial because it forces $\bx$ to be
a vertex.
For the inductive step, assume $\bx$ is a point of $\hcpn$.  If $\bx$ is an
interior point, we are done; if not, $\bx$ lies in some facet of $\hcpn$.
By Lemma \secbc.2, this facet is either a half cube of lower dimension, in 
which case we are done by induction, or it is a simplex of lower dimension,
in which case we appeal to the known corresponding result for simplices.
\qed\enddemo

\head \secc. The topology of the half cube \endhead

In \S\secc, we examine some topological properties of the half cube.  The
half cube itself is not interesting topologically, since it is homeomorphic
to a ball.  However, the geometric
results of \S\secb\ can be used to assemble the $k$-faces of the half cube 
$\hcpn$ into a CW complex $C_n$ in a manner that mirrors its geometric 
properties (Theorem \seccb.2).  The complex $C_n$ has a subcomplex $C_{n, k}$
for each $3 \leq k \leq n$, which is obtained by deleting the
interiors of all half cube shaped faces of dimension at least $k$.
Using a combination of discrete Morse theory and the theory of CW complexes,
we show (Theorem \seccd.2) that the reduced homology of $C_{n, k}$ is
concentrated in degree $k-1$.

\subhead \seccb. The half cube as a CW complex \endsubhead

We now recall the definition of a finite regular CW complex; full details 
may be found in \cite{{\bf 20}, \S8}.

If $X$ and $Y$ are topological spaces with $A \subset X$ and $B \subset Y$,
we define a continuous map $g : (X, A) \ra (Y, B)$ to be a continuous
map $g : X \ra Y$ such that $g(A) \subseteq B$.  If, furthermore,
$g|_{X - A} : X - A \ra Y - B$ is a homeomorphism, we call $G$ a {\it relative
homeomorphism}.

An {\it $n$-cell}, $e = e^n$ is a homeomorphic copy of the open $n$-disk
$D^n - S^{n-1}$, where $D^n$ is the closed unit ball in Euclidean $n$-space
and $S^{n-1}$ is its boundary, the unit $(n-1)$-sphere.  We call $e$ a 
{\it cell} if it is an $n$-cell for some $n$.

If a topological space $X$ is a disjoint union of cells $X = \bigcup \{e :
e \in E\}$, then for each $k \geq 0$, we define the {\it $k$-skeleton}
$X^{(k)}$ of $X$ by $$
X^{(k)} = \bigcup \{e \in E : \dim(e) \leq k\}
.$$

The CW complexes we consider in this paper are all finite, which means that
we can give the following abbreviated definition.

\definition{Definition \seccb.1}
A {\it CW complex} is an ordered triple $(X, E, \Phi)$, where $X$ is a 
Hausdorff space, $E$ is a family of cells in $X$, and $\{\Phi_e : e \in E\}$
is a family of maps, such that 
\item{(i)}{$X = \bigcup \{e : e \in E\}$ is a disjoint union;}
\item{(ii)}{for each $k$-cell $e \in E$, the map $\Phi_e : 
(D^k, S^{k-1}) \ra (e \cup X^{(k-1)}, X^{(k-1)})$ is a relative homeomorphism.}

A {\it subcomplex} of the CW complex $(X, E, \Phi)$ is a triple $(|E'|, E',
\Phi')$, where $E' \subset E$, $$
|E'| := \bigcup \{ e : e \in E' \} \subset X
,$$ $\Phi' = \{\Phi_e : e \in E'\}$, and $\Im\  \Phi_e \subset |E'|$ for every
$e \in E'$.

The {\it Euler characteristic} $\chi(X)$ of a finite CW
complex $(X, E, \Phi)$ is given by $$
\chi(X) = \sum_{i \geq 0} (-1)^i \a_i
,$$ where $\a_i$ is the number of $i$-cells in $E$.  
\enddefinition

The complexes considered here have the 
property that the maps $\Phi_e$ (regarded as mapping to their images) are 
all homeomorphisms.  Such CW complexes are called {\it regular}.

\proclaim{Theorem \seccb.2}
Let $X$ be the half cube $\hcpn$ regarded as a subspace of $\real^n$, and
let $E$ be the union of the following three sets: 
\item{\rm (a)}{the set of vertices of $\hcpn$;}
\item{\rm (b)}{the set of interiors of all $k$-faces of $\hcpn$ for all 
$0 < k < n$;}
\item{\rm (c)}{the interior of $\hcpn$ itself.}

Then $(X, E, \Phi)$ is a regular CW-complex, where the maps $\Phi_e$ are the
natural identifications.
\endproclaim

\demo{Proof}
It is a standard result that the faces of a convex polytope form a regular
CW complex.  The result now follows from Theorem \secbc.8.
\qed\enddemo

\definition{Definition \seccb.3}
Let $E$ be as in Theorem \seccb.2, and let $k \geq 3$.
We define the subset $E(k)$ of $E$ to consist of all elements $e \in E$ 
except those for which $\bar{e}$ is an isometric copy of an $l$-dimensional
half-cube for some $l \geq k$.
\enddefinition

Note that the condition on $\bar{e}$ in the above definition excludes 
precisely those elements $e$ for which $\bar{e} = \Pi(L(\bv, S))$ for some
$\bv \in \Psi^+_n$ and $|S| \geq k$.

\example{Example \seccb.4}
If $k = 3$ in Definition \seccb.3, the effect is to delete the interiors of
all the faces of $\hcpn$ that are isometric to $k$-dimensional half cubes, 
including the tetrahedral faces of this type.  However, the tetrahedral faces
corresponding to $4$-cliques of the form $K(\bv', S)$ are retained, as are
all simplex shaped faces of higher dimensions.  

The case $k = 4$ differs only from the case $k = 3$ in that it includes 
the tetrahedral faces corresponding to $4$-cliques of the form $L(\bv, S)$.
Proposition \secbb.6 then shows that $C_{n, 4}$ is the clique complex
of the half cube graph $\hcgn$.

At the other extreme, if $k > n$ in Definition \seccb.3, we have $E(k) = E$.
If $k = n$, then $E(k)$ includes all but one element of $E$, namely the
$n$-dimensional interior of $\hcpn$.
\endexample

\proclaim{Proposition \seccb.5}
For each $3 \leq k \leq n$, the CW complex $(X, E, \Phi)$ of Theorem
\seccb.2 has a subcomplex $(|E'|, E', \Phi')$, where $E' = E(k)$.  We denote
this subcomplex by $C_{n, k}$.
\endproclaim

\demo{Proof}
Suppose that $\Pi(V)$ and $\Pi(V')$ are faces of $\hcpn$ with respective 
vertex sets $V$ and $V'$.  Suppose in addition that $V \subset V'$.  If
$V$ is of the form $L(\bv, S)$ for some $\bv \in \Psi^+_n$ and $S \subset \bn$,
then $V'$ must also be of the form $L(\bu, S')$, where $|S'| > |S|$, because
the boundary of a simplex shaped face consists of simplex shaped faces of
lower dimension.  It follows that if the interior of $\Pi(V)$ is missing from 
$E'$, then the interior of $\Pi(V')$ is also missing from $E'$, which
implies the statement.
\qed\enddemo

\subhead \seccc. Discrete Morse Theory \endsubhead

Discrete Morse theory, which was introduced by Forman \cite{{\bf 11}}, is a
combinatorial technique for computing the homology of CW complexes.  
By building on work of Chari \cite{{\bf 8}}, Forman later produced a version
of discrete Morse theory based on acyclic matchings in Hasse diagrams
\cite{{\bf 12}}.  This version of the theory plays a key role in computing
the homology of $C_{n, k}$.

\definition{Definition \seccc.1}
Let $K$ be a finite regular CW complex.  
A {\it discrete vector field} on $K$ is
a collection of pairs of cells $(K_1, K_2)$ such that 
\item{(i)}{$K_1$ is a face of $K_2$ of codimension $1$ and}
\item{(ii)}{every cell of $K$ lies in at most one such pair.}

We call a cell of $K$ {\it paired} if it lies in (a unique) one of the above
pairs, and {\it unpaired} otherwise.

If $V$ is a discrete vector field on a regular CW complex $K$, we define
a $V$-path to be a sequence of cells $$
\a_0, \be_0, \a_1, \be_1, \a_2, \ldots, \be_r, \a_{r+1}
$$ such that for each $i = 0, \ldots, r$, (a) each of $\a_i$ and $\a_{i+1}$
is a codimension $1$ face of $\be_i$, 
(b) each $(\a_i, \be_i)$ belongs to $V$ and 
(c) $\a_i \ne \a_{i+1}$ for all $0 \leq i \leq r$.  If
$r \geq 0$, we call the $V$-path {\it nontrivial}, and if $\a_0 = \a_{r+1}$,
we call the $V$-path {\it closed}.  Note that all the faces $\a_i$ have
the same dimension, $p$ say, and all the faces $\be_i$ have dimension $p+1$.
\enddefinition

Let $P$ be the set of cells of $K$, together with the empty cell
$\emptyset$, which we consider to be a cell of dimension $-1$.  The set $P$
becomes a partially ordered set under inclusion.  Let $H$ be the Hasse diagram
of this partial order.  We regard $H$ as a directed graph, in which all edges
point towards cells of larger dimension.

Suppose now that $V$ is a discrete vector field on $K$.  We define $H(V)$ to
be the directed graph obtained from $H$ by reversing the direction of an arrow
if and only if it joins two cells $K_1 \subset K_2$ for which $(K_1, K_2)$
is one of the pairs of $V$.  If the graph $H(V)$ has no directed cycles, we
call $V$ an {\it acyclic (partial) matching} of the Hasse diagram of $K$.

\proclaim{Theorem \seccc.2 (Forman)}
Let $V$ be a discrete vector field on a regular CW complex $K$.
\item{\rm (i)}{There are no nontrivial closed $V$-paths if and only if $V$ is
an acyclic matching of the Hasse diagram of $K$.}
\item{\rm (ii)}{Suppose that $V$ is an acyclic partial matching of the Hasse 
diagram of $K$ in which the empty set is unpaired.  Let $u_p$ denote the number
of unpaired $p$-cells.  Then $K$ is homotopic to a CW complex
with exactly $u_p$ cells of dimension $p$ for each $p \geq 0$.}
\endproclaim

\demo{Proof}
Part (i) is \cite{{\bf 12}, Theorem 6.2} and part (ii) is \cite{{\bf 12}, Theorem 6.3}.
\qed\enddemo

\proclaim{Proposition \seccc.3}
Let $V_k$ be the set of pairs $(\Pi(K(\bv', S)), \Pi(K(\bv', S')))$ 
of cells in $C_{n,k}$ such that $|S| \geq k$, $S \subset S'$ and $S' 
\backslash S = \{n\}$.  Then:
\item{\rm (i)}{$V_k$ is an acyclic matching on the Hasse diagram
of $C_{n, k}$;}
\item{\rm (ii)}{$C_{n, k}$ is homotopic to a CW complex with no cells
in dimension $p$ for any $p \geq k$; in particular, $C_{n, k}$ has
zero homology in dimensions $p \geq k$;}
\item{\rm (iii)}{the homology of $C_{n, k}$ over $\zed$ is concentrated
in degree $k-1$ and is free.}
\endproclaim

\demo{Proof}
We have already seen that $C_{n, k}$ is a regular CW complex.
The hypothesis $|S| \geq k$ together with Lemma \secbb.2 show that $V$ is a 
discrete vector field.  To complete the proof of (i), it suffices by Theorem
\seccc.2 (i) to show that there are no nontrivial closed $V$-paths.

Suppose for a contradiction that $$
\a_0, \be_0, \a_1, \be_1, \a_2, \ldots, \be_r, \a_{r+1}
$$ is such a path, where $\a_{r+1} = \a_0$, $r \geq 0$, and for all 
$0 \leq i \leq r$ we have $\a_i \ne \a_{i+1}$ and $\a_i \subset \be_i$.
By the definition of $V$-path, each of the $(\a_i, \be_i)$
is a pair of $V$.  By definition of $V$, we have $\a_0 = K(\bv', S)$ and
$\be_0 = K(\bv', S \cup \{n\})$, where $\bv' \in \Psi^-_n$ and $n \not\in S$.
Since $\a_1$ is a codimension one face of $\be_0$ that is different from
$\a_0$, we must have $\a_1 = K(\bv', T)$, where $T = S \cup \{n\} \backslash
\{i\}$ for some $1 \leq i < n$.  The fact that $n \in T$ then shows that
$(\a_1, \be_1)$ cannot lie in $V$, which is the contradiction required to 
prove (i).

To prove (ii), we first observe that every cell $\Pi(K(\bv', S))$ for which
$|S| \geq k+1$ is involved in a pair of $V$, whether or not $n$ is an element
of $S$.  It follows that every cell of $C_{n, k}$ of dimension at least $k$
is paired.  Since the empty set (of dimension $-1$) is unpaired in $V$,
the conclusion of (ii) now follows from (i) and Theorem \seccc.2 (ii).

Part (iii) follows from (ii) together with \cite{{\bf 20}, 
Corollary 8.40 (iii)}.
\qed\enddemo

\subhead \seccd. Cellular homology \endsubhead

Cellular homology is a convenient theory for computing homology groups of
CW complexes.  We do not recall the full definition here, but instead refer 
the reader to \cite{{\bf 20}, \S8}.  However, the basic idea is to introduce,
for a CW complex $X$, a free $R$-module $W_k$, where $R$ may be the integers
or a field.
The rank of the $R$-module $W_k$ is the number of $k$-cells in
$X$, and the complex $W_*$ may be equipped with differentials $\pd_*$ so that
the homology of $W_*$ is the (singular) homology, $H_*(X)$, of $X$.  One
may optionally consider $\emptyset$ to be the unique $(-1)$-cell of $X$,
which gives rise to reduced cellular homology.

\proclaim{Lemma \seccd.1}
Let $X$ be a CW complex and let $Y$ be a subcomplex of $X$.  Suppose that
$k \geq 0$ and that for $i \in \{k, k+1\}$, every $i$-cell of $X$ is also
an $i$-cell of $Y$.  Then we have $H_k(X) \cong H_k(Y)$.
\endproclaim

\demo{Proof}
This holds because the chain groups and differentials of the two complexes
agree in degrees $k$ and $k+1$.
\qed\enddemo

\proclaim{Theorem \seccd.2}
For all $3 \leq k \leq n$, the reduced homology of the CW complex $C_{n, k}$ is
concentrated in degree $k-1$, and is free over $\zed$ in degree $k$.
\endproclaim

\demo{Proof}
By Proposition \seccc.3 (iii), the reduced homology of $C_{n, k}$ vanishes
in degrees $k$ and higher and is free in degree $k$.
The complex $C_{n, k}$ agrees with $X$
in degrees $k-1$ and lower, so Lemma \seccd.1 shows that the homology (and,
therefore, the reduced homology) of $C_{n, k}$ agrees with that of $X$ in 
degrees $k-2$ and lower.  Since $X$ is closed, bounded and convex, $X$ is 
contractible and has trivial reduced homology.
\qed\enddemo

\head \secd. Nonvanishing homology \endhead

In \S\secd, we use combinatorial techniques to compute the ranks
of the nonzero homology groups of $C_{n, k}$ (Theorem \secda.2).
These ranks can be readily expressed in terms of generating functions,
or by more explicit formulas (Theorem \secda.5), and using these, we can
verify the coincidence of the $(k-1)$-st Betti number of $C_{n, k}$ with
the $(k-2)$-nd Betti number of a manifold arising in the theory of hyperplane
arrangements (Corollary \secda.6).  We then show that the complex $C_{n, k}$
admits a group of cellular automorphisms arising from the action of the
Coxeter group $G$ of type $D_n$ (Theorem \secdb.3), which endows the 
nonzero homology groups with the structure of $G$-modules.

\subhead \secda. Betti numbers \endsubhead

In \S\secda, we determine the dimensions of the nonzero homology groups
of $C_{n, k}$ appearing in Theorem \seccd.2.  
This is made easy by the following well-known result.

\proclaim{Lemma \secda.1 \cite{{\bf 16}, Theorem 2.44}}
Let $(X, E, \Phi)$ be a finite CW complex.  Then we have $$
\chi(X) = \sum_{i \geq 0} (-1)^i \rank H_i(X)
.$$
\qed\endproclaim

\proclaim{Theorem \secda.2}
The rank of the $(k-1)$-st homology group of $C_{n, k}$ is
given by $$
\rank H_{k-1}(C_{n, k}) = 
\sum_{i = k}^n (-1)^{k+i} 2^{n-i} {n \choose i}
.$$
\endproclaim

\demo{Proof}
Let $\chi(C_{n, k})$ denote the Euler characteristic of the complex
$C_{n, k}$.  From Theorem \seccd.2 and Lemma \secda.1, we see that 
the $(k-1)$-st homology of $C_{n, k}$ satisfies $$
1 + (-1)^{k-1} \rank H_{k-1}(C_{n, k}) = \chi(n, k)
,$$ where the extra ``$1$'' appears because we are considering nonreduced
homology.

By the proof of Theorem \seccd.2, the reduced homology of $C_{n, n+1}$ is 
trivial, which
implies that $\chi(n, n+1) = 1$.  Let $n(i, k)$ be the number of $i$-cells 
in $C_{n, n+1}$ that do not lie in $C_{n, k}$.  If $i < k$, we have 
$n(i, k) = 0$.  On the other hand, if $i \geq k$, we have $$
n(i, k) = 2^{n-i} {n \choose i}
,$$ by Theorem \secbc.6 (v).  By comparing $\chi(n, k)$ with $\chi(n, n+1)$,
we then have $$
\chi(n, k) + \sum_{i = k}^n (-1)^i n(i, k) = 1
.$$

The result now follows by equating the two above expressions for $\chi(n, k)$.
\qed\enddemo

Recall from \cite{{\bf 7}} that the $k$-equal real hyperplane arrangement
$V_{n, k}^\real$ is the set of points $(x_1, \ldots, x_n) \in \real^n$
such that $x_{i_1} = x_{i_2} = \cdots = x_{i_k}$ for some set of indices
$1 \leq i_1 < i_2 < \cdots < i_k \leq n$.  The manifold $M_{n, k}^\real$
is defined to be the complement $\real^n - V_{n, k}^\real$.  
The Betti numbers of $M_{n, k}^\real$ are
special cases of the entries of a particular Pascal-like triangle, $T(n, k)$,
defined as follows.

\definition{Definition \secda.3 \cite{{\bf 21}, sequence A119258}}
Let $n \geq 0$ and $0 \leq k \leq n$ be integers.  We define the sequence
$T(n, k)$ by the conditions
\item{(i)}{$T(n, 0) = T(n, n) = 1$;}
\item{(ii)}{for $0 < k < n$, $T(n, k) = 2 T(n-1, k-1) + T(n-1, k)$.}

We interpret $T(n, k) = 0$ if $k < 0$ or $k > n$, which implies that the 
recurrence relation in (ii) also holds for $k = 0$.
\enddefinition

\example{Example \secda.4}
The entries of $T(n, k)$ for $0 \leq n \leq 6$ and $0 \leq k \leq n$ are
as follows: $$\matrix
 & & & & & & 1 & & & & & & \cr
 & & & & & 1 & & 1 & & & & & \cr
 & & & & 1 & & 3 & & 1 & & & & \cr
 & & & 1 & & 5 & & 7 & & 1 & & & \cr
 & & 1 & & 7 & & 17 & & 15 & & 1 & & \cr
 & 1 & & 9 & & 31 & & 49 & & 31 & & 1 & \cr
1 & & 11 & & 49 & & 111 & & 129 & & 63 & & 1 \cr
\endmatrix $$
\endexample

The next result collects various formulae for sequence $T(n, k)$, and
turns out to be very closely related to Theorem \secda.2.
Some of these formulae are familiar from other contexts.

\proclaim{Theorem \secda.5}
Let $T(n, k)$ be the sequence of Definition \secda.3.
\item{\rm (i)}{The generating function $\sum_{n = 0}^\infty T(n, n-k) x^n$
is given by $$
{x^k \over {(1-2x)^k (1-x)}}
.$$}
\item{\rm (ii)}{We have $$
T(n, k) = \sum_{i = n-k}^n (-1)^{n-k-i} 2^{n-i} {n \choose i}
.$$}
\item{\rm (iii)}{We have $$
T(n, k) = \sum_{i = n-k}^n {n \choose i} {{i-1} \choose {n-k-1}}
,$$ where we interpret ${{-1} \choose {-1}}$ to mean $1$.}
\endproclaim

\demo{Proof}
To prove (i), we note that if $k > 0$, we have the identity $$
(1-2x) {{x^k} \over {(1-2x)^k}} = 
x {{x^{k-1}} \over {(1-2x)^{k-1}}}
.$$  If $n - k > 0$, this expresses the recurrence relation $$
T(n, n-k) - 2 T(n-1, n-1-k) = T(n-1, n-k)
,$$ and this holds by definition of $T(n, n-k)$ since $0 < n-k < n$.  It 
is clear that the coefficient of $x^k$ in $$
{x^k \over {(1-2x)^k (1-x)}}
$$ is $1$, which shows that $T(k, 0) = 1$ for all $k$, as required.  
Finally, if $k = 0$, the coefficient of $x^n$ in $(1-x)^{-1}$  is $1$, 
showing that $T(n, n) = 1$ for all $n$, completing the proof of (i).

Let $T_1(n, k)$ (respectively, $T_2(n, k)$) denote the sequence 
defined by the statement of (ii) (respectively, (iii)).

We prove that $T_1(n, k) = T(n, k)$ by induction on $n$.  The base case,
$n = k = 0$, is trivial, so suppose that $n > 0$.  It is then enough to prove
that $$
T_1(n, k) - T_1(n-1, k) = 2 T_1(n-1, k-1)
.$$  If $L$ is the left hand side of this equation, we have $$\eqalign{
L &= 
\sum_{i = n-k}^n (-1)^{n-k-i} 2^{n-i} {n \choose i}
- \sum_{i = n-1-k}^{n-1} (-1)^{n-1-k-i} 2^{n-1-i} {n-1 \choose i} \cr
&= 
\sum_{i = n-k}^n (-1)^{n-k-i} 2^{n-i} {n \choose i}
- \sum_{i = n-k}^n (-1)^{n-k-i} 2^{n-i} {n-1 \choose i-1} \cr
&= 
\sum_{i = n-k}^n (-1)^{n-k-i} 2^{n-i} {n-1 \choose i} \cr
&= 
2 \sum_{i = (n-1)-(k-1)}^{n-1} (-1)^{(n-1)-(k-1)-i} 2^{n-1-i} 
{n-1 \choose i},\cr
}$$ which equals $2T_1(n-1, k-1)$, proving (ii).

To prove (iii), we show that $T_1(n, k) = T_2(n, k)$ by induction on $k$.
The base case, $k = 0$, follows easily.  For the inductive step, 
we recall that Bj\"orner and Welker \cite{{\bf 7}, 7.5 (ii)} show 
that $$
T_2(n, n-(k+1)) + T_2(n, n-k) = 2^{n-k} {n \choose k}
;$$ they attribute this observation to V. Strehl. 
(Note that this identity also holds for $k = n-1$, because
$T_2(n, 0) = 1$ and $T_2(n, 1) = 2n-1$.)  Equivalently, we may
write $$
T_2(n, k-1) + T_2(n, k) = 2^k {n \choose k}
.$$  The corresponding identity $$
T_1(n, k-1) + T_1(n, k) = 2^k {n \choose k}
$$ is immediate from (ii), and this completes the proof.
\qed\enddemo

\proclaim{Corollary \secda.6}
If $3 \leq k \leq n$, the number $T(n, n-k)$ is equal both to 
\item{\rm (i)}{the $(k-1)$-st Betti number of
the complex $C_{n, k}$ and}
\item{\rm (ii)}{the $(k-2)$-nd Betti number of the complement $M_{n, k}^\real$
of the $k$-equal real hyperplane arrangement.}
\endproclaim

\demo{Proof}
Part (i) follows from comparing Theorem \secda.5 (ii) with Theorem \secda.2.
Part (ii) follows from \cite{{\bf 7}, \S7.5}.
\qed\enddemo

\subhead \secdb. Homology representations \endsubhead

\definition{Definition \secdb.1}
Let $\Aut(\hcgn)$ be the automorphism group of the half cube graph $\hcgn$ 
defined in \S\secbb.

Regard the group $\zed_2 \wr S_n$ as acting as signed permutations on $n$
objects.  Let $W(D_n)$ be the subgroup of $\zed_2 \wr S_n$ consisting
of all elements that effect an even number of sign changes.
\enddefinition

The group $W(D_n)$ is isomorphic to the Coxeter group of type $D_n$; see
\cite{{\bf 17}, \S2.10} for more details.

\proclaim{Lemma \secdb.2}
\item{\rm (i)}
{The group $W(D_n)$ acts as signed permutations on each of the sets
$\Psi^+_n$ and $\Psi^-_n$.  This action induces an embedding of $W(D_n)$
into $\Aut(\hcgn)$.}
\item{\rm (ii)}
{There is an injective group homomorphism $\phi$ from $W(D_n)$ to the
orthogonal group $O_n(\real)$
with the property that whenever $g \in W(D_n)$ and $\psi \in \Psi_n$, 
we have $\phi(g).\psi = g(\psi)$.}
\endproclaim

\demo{Proof}
Part (i) holds because the signed permutation action (a) fixes each of the
two sets $\Psi^\pm_n$ setwise, and (b) respects Hamming distance.

For part (ii), recall from \cite{{\bf 15}, Proposition 5.4} that the two sets
$\Psi^\pm_n$ are the weights of the two spin modules for the simple Lie 
algebra of type $D_n$ over $\complex$.  The action of $W(D_n)$ on these
sets via orthogonal transformations in $\real^n$
is given by \cite{{\bf 15}, Proposition 3.6, Lemma 5.3}.  This action restricts
to the action as signed permutations by the remarks in \cite{{\bf 17}, \S2.10}.
\qed\enddemo

\proclaim{Theorem \secdb.3}
Let $G_n$ be the full subgroup of $O(\real^n)$ stabilizing $\hcpn$ setwise,
and consider $W(D_n)$ as a subgroup of $G_n$ as in Lemma \secdb.2.  
\item{\rm (i)}{The groups $W(D_n)$ and $G_n$ permute the $k$-faces of $\hcpn$ 
for each $k$.}
\item{\rm (ii)}{If $n > 4$, the $G_n$-orbits of $k$-faces coincide with
the $W(D_n)$ orbits.  There are two orbits of $k$-faces 
for $3 \leq k < n$, and one orbit of $k$-faces otherwise.  In the case 
where there are two orbits, the orbits are distinguished by whether they 
correspond to a set $K(\bv', S)$ or a set $L(\bv, S)$ in Theorem \secbc.6.}
\item{\rm (iii)}{If $n = 4$, there is one $G_n$-orbit of $k$-faces for
each $0 \leq k \leq n$. The $W(D_n)$-orbits coincide with the $G_n$-orbits,
except that there are two $W(D_n)$-orbits of $3$-faces; these 
are distinguished by whether they correspond to a set $K(\bv', S)$ or a 
set $L(\bv, S)$ in Theorem \secbc.6.}
\item{\rm (iv)}{The groups $W(D_n)$ act naturally on the nonzero
homology groups $H_{k-1}(C_{n, k})$ for all $3 \leq k \leq n$.}
\endproclaim

\demo{Proof}
Let $g \in W(D_n)$ and consider $g$ acting as a signed permutation.  If
$\bv' \in \Psi^-_n$, $\bv \in \Psi^+_n$ and $S \subseteq \bn$, 
we have $$
g(K(\bv', S)) = K(g.\bv', g(S))
$$ and $$
g(L(\bv, S)) = L(g.\bv, g(S))
,$$ where $g(S)$ refers to $g$ acting on the set $S$ as an ordinary
permutation, ignoring the sign changes.  
Since $g.\bv' \in \Psi^-_n$ and $g.\bv \in 
\Psi^+_n$, $g$ induces a permutation of the subsets of $\Psi_n$ indexing
the faces of $\hcpn$, as in Theorem \secbc.6.  
In particular, we see that the action of $W(D_n)$ both (a) permutes
the $k$-faces of $\hcpn$ for each $k$ and (b) sends faces to faces of the 
same type ($K$-type or $L$-type).

We observe that the vertices of $\hcpn$ are the only points $\bx$ in 
$\hcpn$ for which $\bx \cdot \bx = n$; this property is inherited from 
the hypercube, for which it is obvious.  Since the group $G_n$ acts by
orthogonal transformations, it must therefore permute the vertices of
$\hcpn$.  In turn, $g$ induces a 
permutation of the faces of $\hcpn$ that respects dimensions, and this
completes the proof of (i).

It remains to verify the assertions about the orbits.  Observe that if 
$\Pi(K(\bv', S))$ and $\Pi(K(\bu', S'))$ are two faces of $\hcpn$ for 
which $|S| = |S'|$,
then we can find an element $x \in S_n$ for which $x.S = S'$.  Furthermore,
the fact that $\bv'$ and $\bu'$ lie in $\Psi^-_n$ 
means that the elements $x.\bv'$ and $x.\bu'$
of $\Psi^-_n$ differ by an even number of sign changes.  It follows that
$\Pi(K(\bv', S))$ and $\Pi(K(\bu', S))$ are conjugate under the 
action of $W(D_n)$, and a similar remark applies to faces $\Pi(L(\bv, S))$ 
and $\Pi(L(\bu, S'))$ in the case where $|S| = |S'|$.

To complete the proof, we must show that if $|S| = |T| = r$, then two 
$k$-faces of
the form $\Pi(K(\bv', S))$ and $\Pi(L(\bv, T))$ cannot be conjugate by an 
element $g \in G_n$.  
Since the action of $g$ respects convex hulls, it must induce a
bijection between the vertex sets $K(\bv', S)$ and $L(\bv, T)$.  These
two sets have cardinalities $k+1$ and $2^{k-1}$, respectively, and the 
only integer solution to the equation $k+1 = 2^{k-1}$ for $k \geq 3$ is 
$k = 3$.  However, the assumption that $n \geq 5$ means that the $3$-face
$\Pi(K(\bv', S))$ occurs as a codimension $1$ face of the $4$-face
$K(\bv', S \cup \{i\})$, where $i \not\in S$, whereas the $3$-face
$L(\bv, T)$ is not a codimension $1$ face of any $4$-face.  This shows
that no such $g$ can exist, completing the proof of (ii).

The proof of (iii) follows the same line of argument as the proof of (ii),
except when $k = 3$.  If $k = 3$, it can be checked by hand that the
two $W(D_4)$-orbits of $3$-faces are as claimed.  However, if $g$ is
the orthogonal transformation given by reflection perpendicular to the
vector $(1, 1, 1, 1)$, we find that $g$ stabilizes the set $\Psi^+_n$, and
thus the half cube $\hcpn$.  However, we also have $$
g(L((1, 1, 1, 1), \{1, 2, 3\})) = K((-1, -1, -1, 1), \{1, 2, 3, 4\})
,$$ which proves that there is one $G_4$-orbit of $3$-faces.

To prove (iv), we first observe that any $g \in W(D_n)$ maps the $i$-skeleton
of the CW complex $C_{n, k}$ to itself for all $0 \leq i \leq n$; this
follows from the definition of $C_{n, k}$ together with fact that the 
action of $W(D_n)$ preserves the type of each face.  The filtration of
a CW complex by its sequence of $i$-skeletons is an example of a cellular
filtration (see \cite{{\bf 20}, Theorem 8.38}), and the property of $W(D_n)$ 
just mentioned means that $g$ acts on $C_{n, k}$ by a cellular map, which
is a homeomorphism by invertibility of $g$.  The remarks before 
\cite{{\bf 20}, Theorem 8.38} also show that $g$ induces an automorphism
$g_*$ on the nonzero homology groups of $C_{n, k}$.
\qed\enddemo

\remark{Remark 4.2.4}
The unexpectedly large symmetry group of $\hcpn$ for $n = 4$ occurs because
in this case, $\hcpn$ is the $4$-dimensional hyperoctahedron, also known
as the ``$16$-cell''.  This is a regular polytope, called a ``cross polytope''
in \cite{{\bf 9}} and denoted by $\be_4$.
\endremark

\subhead \secdc. Concluding remarks \endsubhead

We conclude with a number of problems and questions which we believe would
be interesting topics for future research.  One obvious natural question is
the following.

\proclaim{Problem \secdc.1}
Give a conceptual explanation for the numerical coincidence of Betti numbers 
in Corollary \secda.6.
\endproclaim

Barcelo et al \cite{{\bf 2}} introduced a combinatorial homotopy theory for
simplicial complexes and graphs, known as $A$-theory.  Work of Babson et
al \cite{{\bf 1}} and unpublished work of Bj\"orner shows that the nonzero
Betti numbers for $M_{n, k}^\real$ also show up in the context of $A$-theory
as the ranks of $A$-groups.  Furthermore, the coincidence between these 
numbers is fairly well understood, and in the simplest case ($k = 3$),
the problem of computing the Betti numbers reduces 
to the computation of the fundamental group of the 
permutahedron after attaching a $2$-cell to each square; see 
\cite{{\bf 3}, Theorem 5.4} for more details.

\proclaim{Problem \secdc.2}
Find a direct relationship between $A$-theory and the complexes $C_{n, k}$.
\endproclaim

The connection between $A$-theory and hyperplane arrangements gives rise
to another formula for the numbers $T(n, n-k)$ in Theorem \secda.5.  Like
the formula in part (iii) of that result, it consists of a sum of products
of positive numbers.  We did not present it here in order to save space; 
however, Barcelo and Smith \cite{{\bf 4}, \S5} give complete details in the
case $k = 3$, and their arguments adapt easily to the case of general $k$.

\proclaim{Problem \secdc.3}
Find an explicit basis for the nonzero homology groups of $C_{n, k}$.
\endproclaim

The two identities mentioned above that express the Betti numbers as sums
of products of positive integers may give a clue to the nature of such a
basis.  An answer to Problem \secdc.3 would also help with the following
problem.

\proclaim{Problem \secdc.4}
Find the characters (over $\complex$) of the representations of $W(D_n)$ 
afforded by the nonzero homology groups.
\endproclaim

Problem \secdc.4 is likely to be difficult in general, even though the
character theory of the group $W(D_n)$ is well understood \cite{{\bf 14}}.
The symmetric group $S_n$, which is a less complicated group than $W(D_n)$,
acts on the homology of the manifold $M_{n, k}^\real$, but (to the best of
our knowledge) not much is known about the representations that arise in 
this way.  It is possible that progress could be made on Problem \secdc.4
via the next question.

\proclaim{Question \secdc.5}
Does the recurrence relation $$
T(n, k) = 2 T(n-1, k-1) + T(n-1, k)
$$ correspond to a branching rule obtained by restricting the homology
representation arising from $C_{n, k}$ to the subgroup $W(D_{n-1})$?
\endproclaim

Since the $(k-1)$-st homology of $C_{n, k}$ is free over $\zed$, a natural
question is

\proclaim{Question \secdc.6}
Is $C_{n, k}$ homotopic to a wedge of $(k-1)$-spheres?
\endproclaim

It seems likely that this question could be answered in the affirmative using
the homology version of Whitehead's Theorem \cite{{\bf 16}, Proposition 
4C.1}.

\head Acknowledgements \endhead

I thank the referee for reading the paper quickly and carefully, for 
simplifying and correcting some of the arguments, and for pointing out that
the homology of $C_{n, k}$ is in fact free over $\zed$.  I also thank Tim
Penttila for some helpful conversations.

\leftheadtext{} \rightheadtext{}

\end